\tolerance=10000
\magnification=1200
\raggedbottom

\baselineskip=15pt
\parskip=1\jot

\def\sk{\vskip 3\jot}

\def\heading#1{\vskip3\jot{\noindent\bf #1}}
\def\label#1{{\noindent\it #1}}

% Reference format

\def\ref#1;#2;#3;#4;#5.{\item{[#1]} #2,#3,{\it #4},#5.}
\def\refinbook#1;#2;#3;#4;#5;#6.{\item{[#1]} #2, #3, #4, {\it #5},#6.} 
\def\refbook#1;#2;#3;#4.{\item{[#1]} #2,{\it #3},#4.}

{
\pageno=0
\nopagenumbers
\rightline{\tt y.arxiv.tex}
\vskip1in

\centerline{\bf On the Enumeration of Interval Graphs}
\vskip0.5in

\centerline{Joyce C. Yang}
\centerline{\tt jcyang.hmc.edu}
\sk

\centerline{Nicholas Pippenger}
\centerline{\tt njp@math.hmc.edu}
\sk

\centerline{Harvey Mudd College}
\centerline{301 Platt Boulevard}
\centerline{Claremont, CA 91711}
\vskip0.5in

\noindent{\bf Abstract:}
We present upper and lower bounds for the number $i_n$ of interval graphs on $n$ vertices.
Answering a question posed by Hanlon, we show that the ordinary generating function 
$I(x) = \sum_{n\ge 0} i_n\,x^n$ for the number $i_n$ of $n$-vertex interval graphs has radius of convergence zero.
We also show that the exponential generating function $J(x) = \sum_{n\ge 0} i_n\,x^n/n!$ has radius of convergence at least $1/2$.
\vskip0.25in

\leftline{{\bf Keywords:} Graph Theory, Enumeration}
\vskip0.25in

\leftline{{\bf MSC Classification:} 05C30}
\vfill\eject
}

\heading{1, Introduction}

An undirected graph is an {\it interval graph\/} if there is a one-to-one correspondence between its vertices and a set of intervals of the real numbers such that two vertices
are adjacent if and only if their corresponding intervals overlap (that is, have a non-empty intersection).
The concept of an interval graph appears first to have been formulated by Haj\'os [H1] in 1957, who posed the problem of determining whether a given graph is an interval graph.
In 1959, the biologist Benzer [B] independently reformulated the problem in connection with the question of whether the observed overlaps among gene fragments were compatible with the hypothesis that the structure of the gene is linear.
The first intrinsic characterization of interval graphs was given in 1962 by Lekkerkerker and Boland [L], who showed that a graph is an interval graph if and only if it is {\it chordal\/} (that is, every cycle of length at least four has a chord) and {\it anasteroidal\/} (that is, among every three vertices, there is one that is adjacent to every path between the other two).
They also gave a characterization in terms of an infinite family of forbidden induced subgraphs.

In 1982, Hanlon [H] addressed the problem of enumerating interval graphs (that is, determining the number $i_n$ of (isomorphism classes of) interval graphs with $n$ vertices).
He solved this problem by giving a system of equations defining a set of formal power series including the (ordinary) generating function $I(x) = \sum_{n\ge 0} i_n\,x^n$ for interval graphs.
This implicit enumeration allowed him to tabulate $i_n$ for $n$ up to $100$ ($i_{100}$ has $57$ decimal digits), but it did not give an explicit formula for $i_n$, or even allow determination of its asymptotic behavior.
Indeed, Hanlon posed the question of whether the radius of convergence of $I(x)$ is positive (that is, whether $i_n$ is eventually bounded by $C^n$ for some constant $C$).

In Section 2, we shall give a lower bound to $i_n$ that grows ``factorially'', and therefore shows that the radius of convergence of $I(x)$ is zero.
In Section 3, we shall give an upper bound to $i_n$ that shows that the exponential generating function $J(x) = \sum_{n\ge 0} i_n\,x^n/n!$ has a radius of convergence that is at least $1/2$.
\sk

\heading{2. Lower Bound}

In this section, we shall prove the lower bound
$$i_{3k} \ge k!/3^{3k}. \eqno(2.1)$$
Since $i_n$ is non-decreasing, (2.1) implies that the coefficients in the formal power series $I(x)$
are bounded below by those of $(1+x+x^2)\sum_{k\ge 0} k!\, x^{3k}/3^{3k}$.
Since this last power series has radius of convergence zero, so does $I(x)$, thus answering Hanlon's question.

To prove (2.1), we shall associate with each permutation $\pi$ of the set $\{1,\ldots, k\}$
an $3$-colored $(3k)$-vertex interval graph $G_\pi$.
This associate will be one-to-one (that is, the permutation $\pi$ can be recovered from the colored graph $G_\pi$).
Since there at only $i_{3k}\,3^{3k}$ distinct $3$-colored $(3k)$-vertex interval graphs,
we have $i_{3k}\,3^{3k} \ge k!$, which is equivalent to (2.1).

Let $\pi$ be a permutation of $\{1,\ldots,k\}$.
We shall construct $3k$ intervals, with endpoints $1,\ldots,6k$.
First we construct $k$ red intervals, $R_j = [3j-1,3j]$ for $1\le j\le k$.
Then we construct $k$ blue intervals, $B_j = [3k+3j-2,3k+3j-1]$ for $1\le j\le k$.
Finally, we construct $k$ white intervals, $W_j = [3j-2,3k+3\pi(j)]$ for $1\le j\le k$.
Let $G_\pi$ be the interval graph whose colored vertices 
$r_1,\ldots,r_k,b_1,\ldots, b_k,w_1,\ldots,w_k$
correspond to the colored intervals just constructed.
It remains to show that the permutation $\pi$ can be recovered from the colored graph $G_\pi$.
Let us define the {\it red degree\/} $\deg_R(w)$ of a white vertex $w$ to be the number of red vertices adjacent to $w$, and the {\it blue degree\/} $\deg_B(w)$  to be the number of blue vertices adjacent to $w$.
Then in the colored graph $G_\pi$, $w_j$ is the unique white vertex such that 
$\deg_R(w_j)=k+1-j$.
Finally, $\pi(j) = k+1 - \deg_B(w_j)$, which shows that $\pi$ can be recovered from $G_\pi$.

It is clear that (2.1) could be improved slightly, by using only $k-1$ red and $k-1$ blue intervals,
for example, and by using a sharper upper bound to the number of colorings used.
But we have not pursued these improvements, as none of them improve the factor $3$ in the relation between $i_{3k}$ and $k!$.
\sk

\heading{3. Upper Bound}

In this section, we shall prove the upper bound
$$i_n \le (2n-1)!!, \eqno(3.1)$$
where $(2n-1)!! = (2n-1)\cdot(2n-3)\cdots3\cdot1$.
We have $(2n-1)!! = (2n)!/2^n\,n!$ and $(2n-1)!! \le 2^n\,n!$.
The last inequality yields
$$\eqalign{
J(x) 
&= \sum_{n\ge 0} i_n\,x^n/n! \cr
&\le \sum_{b\ge 0} 2^n\,x^n \cr
&= 1/(1-2x), \cr
}$$
which implies, because the coefficients of $J(x)$ are non-negative, that the radius of convergence of $J(x)$ is at least $1/2$.
It remains an open question to determine if $J(x)$ has a larger radius of convergence, or or is perhaps even an entire function.

To prove (3.1), we observe that the $2n$ endpoints of the $n$ intervals in the representation of an $n$-vertex interval graph can, without loss of generality, be taken to be the $2n$ positive integers 
$1, \ldots, 2n$.
The $n$ intervals then correspond to a partition of these $2n$ integers into
$n$ blocks, each containing two integers that are the endpoints of an interval.
This partition can be specified by first choosing the mate of $1$ from among  the $2n-1$ integers $2, \ldots 2n$ (which can be done in $2n-1$ ways), then choosing the mate of the smallest as-yet-unpaired integer from among the $2n-3$ greater as-yet-unpaired integers (which can be done in $2n-1$ ways), and continuing in this way until all $2n$ integers have been partitioned into $n$ mated pairs.
This can thus be done in $(2n-1)\cdot(2n-3)\cdots3\cdot1 = (2n-1)$ ways.
Since every interval graph corresponds to at least one partition, we have established (3.1).
\sk

\heading{4. Conclusion}

In view of Stirling's formula,
one way of roughly stating our results is
$${1\over 3}n\log n + O(n) \le \log i_n \le n\log n + O(n).$$
It remains an open problem to bring the coefficients $1/3$ and $1$ in the lower and upper bounds closer together, perhaps even to obtain an estimate of the form
$$\log i_n = C\,n\log n + O(n)$$
for some constant $C$.
\sk

\heading{5. References}

\ref B; S. Benzer;
``On the Topology of the Genetic Fine Structure'';
Proc.\ Nat.\ Acad.\ Sci.\ USA; 45:11 (1959) 1607.

\ref H1; G. Haj\'os;
``\"Uber eine Art von Graphen'';
Int.\ Math.\ Nachr.; 11 (1957) 65.

\ref H2; Ph.\ Hanlon;
``Counting Interval Graphs'';
Trans.\ Amer.\ Math.\ Soc.; 272:21 (1982) 383--426.

\ref L; C. G. Lekkerkerker and J. Ch.\ Boland;
``Representation of a Finite Graph by a Set of Intervals on the Real Line'';
Fund.\ Math.; 51 (1962) 45--64.

\bye